\documentclass{amsart} 
\usepackage{amssymb, amsmath}%%% option packages
\usepackage[dvips]{graphicx}%%% option packages
\usepackage{amsfonts}
\usepackage{latexsym}
\usepackage{color}

\newtheorem{theorem}{Theorem}	
\newtheorem{lemma}{Lemma}%[section]		

\newtheorem{definition}{Definition}	
\newtheorem{example}{Example}

\title[Kato's chaos created by quadratic mappings]
{
Kato's chaos created by quadratic mappings 
associated with spherical orthotomic curves
}
\author{Takashi Nishimura%\thanks{Corresponding author.} 
}
\address{
Takashi Nishimura: Research Institute of Environment 
and Information Sciences,  
Yokohama National University, 
Yokohama 240-8501, Japan}
\email{nishimura-takashi-yx@ynu.jp}
\begin{document}
\begin{abstract}
{
Singular quadratic mappings creating 
Kato's chaos are given.   
} 
\end{abstract}
%%%%%%%%%%%%%%%%%%%%%%%%%%%%%%%%%% 
\maketitle  
\section{Introduction}\label{section 1} 
\noindent
Throughout this paper, let $n$ be a non-negative integer. 
Moreover, let $S^n, D^{n+1}$ be the unit sphere and 
the unit disk of $\mathbb{R}^{n+1}$ respectively.   
For any point $P\in S^n$, define the quadratic mapping 
$\Phi_P: \mathbb{R}^{n+1}\to \mathbb{R}^{n+1}$ as follows:  
\[
\Phi_P(x)=2\left(x\cdot P\right)x-P, 
\]     
where the dot in the center stands for the standard scalar product 
of two vectors $x$ and $P$ of $\mathbb{R}^{n+1}$.   
%\par 
The mapping $\Phi_P$, which is the subject of our research in this paper, 
can be naturally obtained by considering 
the relation between spherical pedal curves and 
spherical orthotomic curves  (see Section \ref{section 4}).
%It is not difficult to prove the following (for the proof 
%of Lemma \ref{lemma 1}, see Section \ref{section 2}).   
\begin{lemma}\label{lemma 1} 
For any $P\in S^n$, the following hold:  
\begin{enumerate}
\item[(1)] $\Phi_P(S^n)\subset S^n$ for any $n\ge 0$. 
\item[(2)] $\Phi_P(S^n)\supset  S^n$ for any $n\ge 1$.      
\item[(3)] $\Phi_P(D^{n+1})=D^{n+1}$ for any $n\ge 0$.   
\end{enumerate}
\end{lemma}
\noindent 
For the proof of Lemma \ref{lemma 1}, 
see Section \ref{section 2}.    
From Lemma \ref{lemma 1}, discrete dynamics for 
$\Phi_P|_{S^n}$ $(n\ge 1)$ and $\Phi_P|_{D^{n+1}}$ $(n\ge 0)$ seems to be  significant to be investigated.   
%\par 
%\medskip 
%We would like to investigate chaotic behavior of iteration of $\Phi_P$.    
\begin{example}
{\rm 
Suppose that $n=1$ and $P=(1,0)$.    
Then, $\Phi_P(x)=(2x_1^2-1, 2x_1x_2)$, where $x=(x_1, x_2)$.   
If $x$ belongs to $S^1$, $x$ may be wirtten as 
$x=\left(\cos\theta, \sin\theta\right)$.    
Then, 
\[
\Phi_P|_{S^1}(\cos\theta, \sin\theta) 
 =  \left(2\cos^2\theta-1, 2\cos\theta\sin\theta\right)   
 =  
\left(\cos 2\theta, \sin 2\theta\right).    
\]
Thus, the restricted mapping $\Phi_P|_{S^n}$ in this case 
is exactly the same mapping given in Chapter 1 
Example 3.4 of Devaney's well-known book \cite{devaney}.  
}
\end{example}
\begin{example}       
{\rm 
Suppose that $n=0$.    Then, $P$ is $1$ or $-1$, and 
$\Phi_P(x)=2x^2-1$ or $-2x^2+1$.   
Define the affine transformation $h_P: \mathbb{R}\to \mathbb{R}$ 
as follows:   
\[
h_P(x)=
\left\{
\begin{array}{rl}
-2x+1 & \quad (\mbox{if }\; P=1) \\ 
2x-1 & \quad (\mbox{if }\; P=-1). 
\end{array}
\right.
\]
Then, in each case, it is easily seen that 
$h_P^{-1}\circ\Phi_P\circ h_P(x)=4x(1-x)$.   
Therefore, in each case, 
$\Phi_P|_{D^1}$ has the same dynamics 
as Chapter 1 Example 8.9 of \cite{devaney}.     
}
\end{example}  
From Examples 1 and 2, it seems meaningful   
to study the chaotic behavior  
of iteration for $\Phi_P|_{S^n}:S^n\to S^n$ $(n\ge 1)$ or 
$\Phi_P|_{D^{n+1}}:D^{n+1}\to D^{n+1}$ $(n\ge 0)$, which is the main purpose of this paper.     
%There are many different definitions of chaos.   
%It seems that the following three are relatively common definitions.  
\begin{definition}\label{definition 1}
{\rm Let $(X, d)$ be a %compact 
metric space with metric $d$ and let 
$f: X\to X$ be a continuous mapping.   
\begin{enumerate}
\item[(1)] 
The mapping $f$ is said to be 
{\it sensitive} 
if there is a positive number $\lambda>0$ such that for any  
$x\in X$ and any neighbourhood $U$ of $x$ in $X$, there exists 
a point $y\in U$ and a non-negative integer $k\ge 0$ such that 
$d(f^k(x), f^k(y)) > \lambda$, 
where 
$f^k$ stands for $\underbrace{f\circ \cdots \circ f}_{\mbox{$k$-tuples}}$.   
\item[(2)] 
The mapping $f$ is said to be {\it transitive} 
if for any non-empty open subsets $U,V\subset X$, there exists a 
positive integer $k> 0$ such that $f^k(U)\cap V\ne \emptyset$.   
\item[(3)] 
The mapping $f$ is said to be {\it accessible} if for any 
$\lambda>0$ and any non-empty open subsets $U,V\subset X$, 
there exist  two points $u\in U$, $v\in V$ and  a 
positive integer $k> 0$ such that $d(f^k(u), f^k(v))\le\lambda$.  
\item[(4)] 
The mapping $f$ is said to be {\it topologically mixing} 
if for any non-empty open subsets $U,V\subset X$, there exists a 
positive integer $k> 0$ such that $f^m(U)\cap V\ne \emptyset$ for any 
$m\ge k$.  
%\item[(5)] 
%The mapping $f$ is said to be {\it chaotic in the sense of Ruelle-Takens }
%(\cite{ruelletakens}) if 
%$f$ is sensitive and transitive.     
\item[(5)] 
The mapping $f$ is said to be {\it chaotic in the sense of Devaney }
(\cite{devaney}) if 
$f$ is sensitive, transitive and the set consisting of periodic points of 
$f$ is dense in $X$. 
\item[(6)] 
The mapping $f$ is said to be {\it chaotic in the sense of Kato }
(\cite{kato}) if 
$f$ is sensitive and accessible.  
\end{enumerate} 
}
\end{definition}  
By definition, it is clear that if a mapping $f: X\to X$ is topologically 
mixing, then it is transitive.   
%; and that if a mapping $f: X\to X$ is chaotic 
%in the sense of Devaney, then it is chaotic in the sense of Ruelle-Takens.     
Moreover, by \cite{kato}, it is known that 
if a mapping $f: X\to X$ is topologically mixing, 
then it is chaotic in the sense of Kato.   
Although Kato's chaos has been well-investigated 
(for instance, see \cite{kato, liwangzhao, wangliangchu}),  
elementary examples  
which is singular and not transitive seems to have been  desired.    
Theorem \ref{theorem 1} gives such examples.  
%elementary concrete examples 
%which is not topologically mixing in any positive dimensions.  
%Moreover, these examples are not transitive in higher dimensions.    
\begin{theorem}\label{theorem 1}
\begin{enumerate}
\item[(1)] 
Let $P$ be a point of $S^1$.  
Then, $\Phi_P|_{S^1}: S^1\to S^1$ is chaotic in the sense of Devaney. 
Moreover, it is chaotic in the sense of Kato. 
\item[(2)] 
Let $P$ be a point of $S^0$.  
Then, $\Phi_P|_{D^1}: D^1\to D^1$ is chaotic in the sense of Devaney. 
Moreover, it is chaotic in the sense of Kato. 
\item[(3)] 
Let $m$ be an integer such that $m\ge 2$.   
Moreover, let $P$ be a point of $S^m$.  
Then, both  $\Phi_P|_{D^m}: D^m\to D^m$  
and $\Phi_P|_{S^m}: S^m\to S^m$ are chaotic in the sense 
of Kato.    
\item[(4)]
Let $m$ be an integer such that $m\ge 2$.   
Moreover, let $P$ be a point of $S^m$.  
Then, 
neigher $\Phi_P|_{D^m}: D^m\to D^m$  
nor $\Phi_P|_{S^m}: S^m\to S^m$ is transitive.   
In particular,  neigher $\Phi_P|_{D^m}: D^m\to D^m$  
nor $\Phi_P|_{S^m}: S^m\to S^m$ is chaotic in the sense of 
Devaney.   
\end{enumerate}
\end{theorem}
\par 
%\noindent 
\bigskip 
%\noindent 
%%%%%%%%%%%%%%%%%%%%%%%%%%%%%%  
This paper is organized as follows.   
In Section \ref{section 2}, 
the proof of Lemma \ref{lemma 1} is given.   
Theorem \ref{theorem 1} is proved in Section \ref{section 3}.   
Section 4 is the appendix where 
how to discover the quadratic mapping $\Phi_P$ is explained.   
%%%%%%%%%%%%%%%%%%%%%%%%%%%%%%%%%%%%%%%%%%%%%%%%%%%%%%%%%%   
%%%%%%%%%%%%%%%%%%%%%%%%%%%%%%%%%%%%%%%%%%%%%%%%%%%%%%%%%% 
\section{Proof of Lemma \ref{lemma 1}}\label{section 2}
%%%%%%%%%%%%%%%%%%%%%%%%%%%%%%%%%%%%%%%%%%%%%%%%%%%%%%%%%% 
\subsection{Proof of the assertion (1) of Lemma \ref{lemma 1}}     
%%%%%%%%%%%%%%%%%%%%%%%%%%%%%%%%%%%%%   
Let $x$ be a point of $S^n$.    Then, $x\cdot x=1$ and we have the following:   
%It sufficies to show that $\Phi_P(x)\cdot \Phi_P(x)=1$ 
%under the assumption 
%$x\cdot x=1$.   
\begin{eqnarray*}
\Phi_P(x)\cdot \Phi_P(x) & = & 
\left(2(x\cdot P)x-P\right)\cdot \left(2(x\cdot P)x-P\right) \\ 
{ } & = & 
4(x\cdot P)^2(x\cdot x)-4(x\cdot P)^2+(P\cdot P) \\ 
{ } & = & 
4(x\cdot P)^2-4(x\cdot P)^2+1=1. \\ 
\end{eqnarray*} 
This completes the proof of the assertion (1).  
\hfill $\Box$  
%\par 
%%%%%%%%%%%%%%%%%%%%%%%%%%%%%% 
\subsection{Proof of the assertion (2) of Lemma \ref{lemma 1}}
%%%%%%%%%%%%%%%%%%%%%%%%%%%%%%     
Let $y$ be a point of $S^n$.    
Suppose that $y\ne -P$.   
Set 
\[
x=\frac{\frac{y+P}{2}}{||\frac{y+P}{2}||}.   
\]
Then, it follows 
\begin{eqnarray*}
2(x\cdot P)x-P 
& = & 
2\left(\frac{\frac{y+P}{2}}{||\frac{y+P}{2}||}\cdot P\right)
\frac{\frac{y+P}{2}}{||\frac{y+P}{2}||}-P \\ 
{ } & = & 
\frac{2}{|| y+P||^2}\left((y\cdot P)+1\right)(y+P)-P \\ 
{ } & = & 
\frac{1}{\left(1+(y\cdot P)\right)}\left((y\cdot P)+1\right)(y+P)-P \\ 
{ } & = & 
(y+P)-P=y.   
\end{eqnarray*}
Next, suppose that $y=-P$.   
Let $x$ be a point of $S^n$ such that $x\cdot P=0$.    
Then, $2(x\cdot P)x-P=-P=y$.    
Therefore, we have the assertion (2).  
\hfill $\Box$     
%\par 
%%%%%%%%%%%%%%%%%%%%%%%%%%%%%%%%%%%%%  
\subsection{Proof of the assertion (3) of Lemma \ref{lemma 1}}     
%%%%%%%%%%%%%%%%%%%%%%%%%%%%%%%%%%%%% 
Let $x$ be a point of $\mathbb{R}^{n+1}$ such that 
$x\cdot x< 1$.    Then, we have 
\[
\Phi_P(x)\cdot \Phi_P(x)  <  
4(x\cdot P)^2-4(x\cdot P)^2+1  =  1.  
\]
Conversely, let $y$ be a point satisfying  
$y\cdot y<1$.   
Notice that in this case 
$(y\cdot P)+1\ge -||y||+1>0$ and 
$1+||y||^2+2(y\cdot P)\ge 1+||y||^2-2||y||=(1-||y||)^2>0$.   
Set 
\[
a=\sqrt{\frac{1+||y||^2+2(y\cdot P)}{2(y\cdot P)+2}}
\mbox{ and }
x=a\frac{\frac{y+P}{2}}{||\frac{y+P}{2}||}.  
\]
Then, 
\begin{eqnarray*}
2(x\cdot P)x-P 
& = & 
2\left(a\frac{\frac{y+P}{2}}{||\frac{y+P}{2}||}\cdot P\right)
a\frac{\frac{y+P}{2}}{||\frac{y+P}{2}||}-P \\ 
{ } & = & 
\frac{2a^2}{|| y+P||^2}\left((y\cdot P)+1\right)(y+P)-P \\ 
{ } & = & 
\frac{2a^2}{\left(1+||y||^2+2(y\cdot P)\right)}
\left((y\cdot P)+1\right)(y+P)-P \\ 
{ } & = & 
(y+P)-P=y.   
\end{eqnarray*}
Therefore, the assertion (3) holds.  
\hfill $\Box$  
%%%%%%%%%%%%%%%%%%%%%%%%%%%%%%%%%%%%%%%%%%%%%%%%%%%%%%%%%%   
%%%%%%%%%%%%%%%%%%%%%%%%%%%%%%%%%%%%%%%%%%%%%%%%%%%%%%%%%%
\section{Proof of Theorem \ref{theorem 1}}\label{section 3}
%%%%%%%%%%%%%%%%%%%%%%%%%%%%%%%%%%%%%%%%%%%%%%%%%%%%%%%%%%
\subsection{Proof of the assertion (1) of Theorem \ref{theorem 1}}
\label{subsection 3.1}     
%%%%%%%%%%%%%%%%%%%%%%%%%%%%%%%%%%%%%
Let $x$ be a point of $S^1$.   Set 
\[
P=(\cos\alpha, \sin\alpha) 
\mbox{ and }
x=(\cos\theta, \sin\theta).  
\]
Then, it is easily seen that  
\begin{eqnarray*}
{ } & { } & \Phi_P(\cos\theta, \sin\theta) \\ 
{ } & = & 
2\left((\cos\alpha, \sin\alpha)\cdot (\cos \theta, \sin\theta)\right)
(\cos\theta, \sin\theta)-
(\cos \alpha, \sin\alpha) \\ 
{ } & = &  
\left(\cos(2\theta-\alpha), \sin(2\theta-\alpha)\right).    
%\leqno{(*)}
\end{eqnarray*}
It follows  
$\Phi_P^k\left(\cos(\theta+\alpha), \sin(\theta+\alpha)\right)=
\left(\cos(2^k\theta+\alpha), \sin(2^k\theta+\alpha)\right)$ 
and therefore,  
by the same argument as in Example 8.6 of \cite{devaney}, 
$\Phi_P$ is chaotic in the sense of Devaney.    
In order to show that $\Phi_P$ is chaotic in the sense of Kato, 
it is sufficient to show that $\Phi_P$ is accessible, which is 
easily seen by the above formula.  
%\[
%\Phi_P(\cos\theta, \sin\theta)= 
%\left(\cos(2\theta-\alpha), \sin(2\theta-\alpha)\right). 
%
%\]   
\hfill 
$\Box$
%%%%%%%%%%%%%%%%%%%%%%%%%%%%%%%%%%%%%
\subsection{Proof of the assertion (2) of Theorem \ref{theorem 1}}
\label{subsection 3.2}     
%%%%%%%%%%%%%%%%%%%%%%%%%%%%%%%%%%%%% 
By Subsection \ref{subsection 3.1} and Example 8.9 of \cite{devaney}, 
$\Phi_P$ is chaotic in the sense of Devaney.    
Moreover, it is easily seen that 
the property of accessibility is preserved by semi-conjugacy.   
Thus, $\Phi_P$ is chaotic in the sense of Kato as well.   
\hfill 
$\Box$
%%%%%%%%%%%%%%%%%%%%%%%%%%%%%%%%%%%%%
\subsection{Proof of the assertion (3) of Theorem \ref{theorem 1}}   
\label{subsection 3.3}  
%%%%%%%%%%%%%%%%%%%%%%%%%%%%%%%%%%%%%
Let ${Q}$ be a point of $S^m$ satisfying $P\ne Q$.    
Set 
\[
{P_Q^{\perp}}=\frac{{Q}-(P\cdot {Q})P}
{||{Q}-(P\cdot Q)P||}.   
\]
Then, it follows $P_Q^{\perp}\in S^m$ and $P\cdot P_Q^{\perp}=0$.   
Let $x$ be a point of the circle 
$S^m\cap (\mathbb{R}P+\mathbb{R}P_Q^{\perp})$.    
%by $S^1(P,Q)$.    
%Let $x$ be a point of  $S^1(P,Q)$.   
Then, $x$ may be written as 
$x=\cos \theta\; P+\sin\theta\;  P_Q^{\perp}$.    
Then, it is easily seen that 
\[
\Phi_P(\cos \theta\; P+\sin\theta\;  P_Q^{\perp})= 
\cos 2\theta\; P+\sin2\theta\;  P_Q^{\perp}.  
\]
Hence, for any non-empty open neighborhood $U$ of $Q$ in $S^n$ 
there exists a positive integer $i$ such that 
the circle $S^m\cap (\mathbb{R}P+\mathbb{R}P_Q^{\perp})$ 
is contained in $\Phi_P^i(U)$.    
Therefore, $\Phi_P|_{S^m}$ is sensitive.   
\par 
Next,  take another point $R$.    
By the same argument as above, it is seen that 
for any non-empty open neighborhood $V$ of $R$ in $S^n$ 
there exists a positive integer $j$ such that 
the circle $S^m\cap (\mathbb{R}P+\mathbb{R}P_Q^{\perp})$ 
is contained in $\Phi_P^j(V)$.  
Set  $k=\max(i,j)$.     
Then, it follows 
\[
P\in \Phi_P^k(U)\cap \Phi_P^k(V).   
\]
Hence, $\Phi_P|_{S^m}$ is accessible.   
\par 
Moreover, 
since $\Phi_P|_{S^m}$ and $\Phi_P|_{D^m}$ are semi-conjugate, 
$\Phi_P|_{D^m}$ is also sensitive and accessible.  
Therefore, both $\Phi_P|_{S^m}$ and $\Phi_P|_{D^m}$ are 
chaotic in the sense of Kato.   
\hfill 
$\Box$
%%%%%%%%%%%%%%%%%%%%%%%%%%%%%%%%%%%%
\subsection{Proof of the assertion (4) of Theorem \ref{theorem 1}} 
\label{subsection 3.4}    
%%%%%%%%%%%%%%%%%%%%%%%%%%%%%%%%%%%%%
%It is sufficient to show that neigher $\Phi_P|_{S^m}$ and  
%$\Phi_P|_{D^m}$ are not topologically mixing.   
Let $Q, R$ be points of $S^m$ so that 
$P, Q, R$ are linear independent.   
Then, $R$ does not belong to the circle 
$S^m\cap (\mathbb{R}P+\mathbb{R}P_Q^{\perp})$ 
where $P_Q^{\perp}$ is the point constructed in 
Subsection \ref{subsection 3.3}.    
Thus, by the argument given in Subsection \ref{subsection 3.3}, 
there exist sufficiently small neighborhoods $U$ (resp., $V$) of 
$Q$ (resp., $R$) in $S^m$ such that 
$\Phi_P^\ell(U)\cap V=\emptyset$ for any $\ell\ge 0$.   
Hence, $\Phi_P|_{S^m}$ is never transitive.   
Since $\Phi_P|_{S^m}$ and $\Phi_P|_{D^m}$ are semi-conjigate, 
even $\Phi_P|_{D^m}$ is not transitive.  
\hfill 
$\Box$
%%%%%%%%%%%%%%%%%%%%%%%%%%%%%%%%%%%%%%%%%%%%%%%%%  
%%%%%%%%%%%%%%%%%%%%%%%%%%%%%%%%%%%%%%%%%%%%%%%%% 
\section{How to discover $\Phi_P$}
\label{section 4}
%%%%%%%%%%%%%%%%%%%%%%%%%%%%%%%%%%%%%%%%%%%%%%%%% 
\subsection{In the plane $\mathbb{R}^2$}
%%%%%%%%%%%%%%%%%%%%%%%%%%%%%%%%%%%%%%%%%%%%%%%%%
In \cite{brucegiblin}, for a given plane unit-speed curve 
$\gamma: I\to \mathbb{R}^2$ and a given point $P\in \mathbb{R}^2$, 
the pedal curve $ped_{\gamma, P}: I\to \mathbb{R}^2$ 
and the orthotomic curve $ort_{\gamma, P}: I\to \mathbb{R}^2$ 
are defined as follows:   
\begin{eqnarray*}
ped_{\gamma, P}(s) & = & P+\left((\gamma-P)\cdot N(s)\right)N(s) \\ 
ort_{\gamma, P}(s) & = & P+2\left((\gamma-P)\cdot N(s)\right)N(s).   
\end{eqnarray*}
Here, $N(s)$ is the unit normal vector to $\gamma$ at $\gamma(s)$.    
Thus, it follows  
\[\frac{ort_{\gamma, P}(s)+P}{2}=ped_{\gamma, P}(s)
\] and thus 
$ort_{\gamma, P}(s)=2ped_{\gamma, P}(s)-P$.   
Therefore, by using the simple mapping 
$F_P: \mathbb{R}^2\to \mathbb{R}^2$ 
defined by 
\[
F_P(x)=2x-P, 
\] 
we have the following:   
\[
ort_{\gamma, P}(s)=F_P\circ ped_{\gamma, P}(s).   
\leqno{(*)}
\]
%%%%%%%%%%%%%%%%%%%%%%%%%%%%%%%%%%%%%%%%%%%% 
\subsection{In the unit sphere $S^n$ %of $\mathbb{R}^{n+1}$ 
}
%%%%%%%%%%%%%%%%%%%%%%%%%%%%%%%%%%%%%%%%%%%% 
In the unit sphere of $\mathbb{R}^{n+1}$, it is desired to have 
an equality like (*) in the plane case.    
In order to do so, for a generic 
unit-speed curve $\gamma: I\to S^n$ and a generic point 
$P\in S^n$, 
%first of all, 
%the unit-speed curve $\gamma: I\to S^n$ need to be defined.   
the pedal curve $ped_{\gamma, P}: I\to S^n$ 
and the orthotomic curve 
$ort_{\gamma, P}: I\to S^n$ need to be defined reasonably.   
In \cite{geomdedicata, demonstratio}, a reasonable definition of unit speed 
curve is given; and then for a unit speed curve $\gamma: I\to S^n$ 
and a generic point $P\in S^n$, $ped_{\gamma, P}: I\to S^n$ is defined 
reasonably.     
Notice that well-definedness of $ped_{\gamma, P}: I\to S^n$ implies  $P\cdot ped_{\gamma, P}(s)\ne 0$ for any $s\in I$ (see \cite{geomdedicata, demonstratio}).    
Thus, we have the relation between $ped_{\gamma, P}: I\to S^n$ 
and $ort_{\gamma, P}: I\to S^n$ as follows. 
\[
\frac{ort_{\gamma, P}(s)+P}{2}= 
\left(P\cdot ped_{\gamma, P}(s)\right)ped_{\gamma, P}(s).   
\]
Hence, we have 
$ort_{\gamma, P}(s)=
2\left(P\cdot ped_{\gamma, P}(s)\right)ped_{\gamma, P}(s)-P$.  
Therefore, both the mapping   
\[
\Phi_P(x)=2(P\cdot x)x-P
\] 
%defined by 
%$\Phi_P(x)=2(P\cdot x)x-P$ and 
and the equality 
\[
ort_{\gamma, P}(s)=\Phi_P\circ ped_{\gamma, P}(s).  
\]
are naturally obtained.   
%%%%%%%%%%%%%%%%%%%%%%%%%%%%%%%%%%%%%%%%%%%%%%%%% 
%%%%%%%%%%%%%%%%%%%%%%%%%%%%%%%%%%%%%%%%%%%%%%%%%%%%
\section*{Acknowledgement}
%\thanks
The author is partially supported 
by JSPS KAKENHI Grant Number 17K05245.   
%%%%%%%%%%%%%%%%%%%%%%%%%%%%%%%%%%%%%%%%%%%%%%%%%%%% 

\end{document}